\documentclass{elsarticle}

\usepackage{amsmath, amsfonts,amsthm,times,graphics}

\begin{document}

\textheight212mm

\begin{frontmatter}

\title{A search for primes $p$ such that
 Euler number $E_{p-3}$ is divisible by $p$} 

\author{Romeo Me\v{s}trovi\'{c}}

\address{Department of Mathematics, Maritime Faculty Kotor,
University of Montenegro\\
Dobrota 36, 85330 Kotor, Montenegro\\
e-mail: romeo@ac.me}

\begin{abstract}
Let $p>3$ be a prime. Euler numbers $E_{p-3}$
first appeared in H. S. Vandiver's work (1940) in connection 
with the first case of Fermat Last Theorem. 
Vandiver proved that $x^p+y^p=z^p$ has no solution for 
integers $x,y,z$ with $\gcd(xyz,p)=1$ if $E_{p-3}\equiv 0\,(\bmod\,p)$.
Numerous combinatorial  congruences
recently obtained by Z.-W. Sun 
and by Z.-H. Sun  
involve the Euler numbers $E_{p-3}$. 
This gives a new  significance to the primes $p$ for which 
$E_{p-3}\equiv 0\,(\bmod\,p)$.

For the computation of residues of Euler numbers $E_{p-3}$ 
modulo a prime $p$, we use the  
congruence which runs significantly faster than
other known congruences involving $E_{p-3}$.
Applying this congruence, 
a  computation via {\tt Mathematica 8} shows that only  
three primes less than $10^7$ 
 satisfy the condition $E_{p-3}\equiv 0\,(\bmod\,p)$
(such primes are 149, 241 and 2946901, and 
they are given as a Sloane's sequence A198245). 
By using related computational results and statistical considerations similar 
to those on search  for Wieferich and Fibonacci-Wieferich 
and Wolstenholme primes, we conjecture that there are infinitely many 
primes $p$ such that $E_{p-3}\equiv 0\,(\bmod\,p)$.
Moreover, we propose a conjecture on the asymptotic estimate of 
number of primes $p$ in an interval $[x,y]$ such 
that $E_{p-3}\equiv A\,(\bmod\,p)$ for some integer $A$
with $|A|\in [K,L]$. 
  \end{abstract}
\end{frontmatter}

{\renewcommand{\thefootnote}{}

\section{Introduction}
{\it Euler numbers}  $E_n$ $(n=0,1,2,\ldots)$ 
(e.g., see \cite[pp. 202--203]{r})
  are integers defined 
recursively by
   $$
 E_0=1, \quad \mathrm{and}\quad\sum_{0\le k\le n\atop k\,\, even}
{n\choose k}E_{n-k} \quad \mathrm{for}\quad n=1,2,3,\ldots
   $$
(it is well known  that $E_{2n-1}=0$ for each $n=1,2,\ldots$).
The first few Euler numbers are $E_0=1, E_2=-1, E_4=5,E_6=-61,E_8=1385,
E_{10}=-50521,E_{12}=2702765,E_{14}=-199360981,E_{16}=19391512145$.
It is well known that Euler numbers can also be defined by the generating 
function
  $$
\frac{2}{e^x+e^{-x}}=\sum_{n=0}^{\infty}E_n\frac{x^n}{n!}.
 $$
It is well known that $E_n=E_n(0)$ ($n=0,1,\ldots$),
where $E_n(x)$ is the classical  Euler polynomial
(see e.g., \cite[p. 61 {\it et seq.}]{sc}).

Recall that {\it Bernoulli numbers} $B_n$ $(n=0,1,2,\ldots)$
are rational numbers defined by the formal identity 
  $$
  \frac{x}{e^x-1}=\sum_{n=0}^{\infty}B_n\frac{x^n}{n!}.
 $$
It is easy to see that $B_n=0$ for odd $n\ge 3$, 
and the first few nonzero terms of $(B_n)$  are $B_0=1$, $B_1=-1/2$, 
$B_2=1/6$, $B_4=-1/30$,  $B_6=1/42$ and $B_8=-1/30$.
It is well known that $B_n=B_n(0)$, where $B_n(x)$ 
is the classical Bernoulli polynomial
(see e.g., \cite[p. 61 {\it et seq.}]{sc}). 

A significance of Euler numbers, and 
especially of  $E_{p-3}$ with a prime  $p$, is closely related to
 Fermat Last Theorem 
(see \cite[Lecture X, Section 2]{r}). In 1850 Kummer (see e.g., 
\cite[Theorem (3A), p. 86 and Theorems (2A)--(2F), pp. 99--103]{r} proved 
that Fermat Last Theorem holds for each {\it regular prime},
that is, for each prime  $p$ that does not divide the numerator of any 
Bernoulli number $B_{2n}$ with $n=1,2,\ldots, (p-3)/2$. 
In 1940 H. S. Vandiver \cite{va} likewise proved for Euler-regular primes.
Paralleling the previous definition of a (irr)regular prime 
(with respect to the Bernoulli numbers)
following Vandiver \cite{va},  a prime $p$ is said to be 
{\it Euler-irregular primes} (shortly $E$-{\it irregular}) if and only if it 
divides at least one of the Euler numbers $E_{2n}$ with $1\le n\le (p-3)/2$. 
Otherwise, that is if $p$ does not divide $E_2,E_4,\ldots,E_{p-3}$, a prime 
$p$ is called $E$-{\it regular}.
The smallest $E$-irregular prime is $p=19$, which divides 
$E(10)=-50521$. The first few $E$-irregular primes are $19,31,43,47,61,67,71,
79, 101,137,139,149,193,223,241$ (with $p=241$ dividing both $E_{210}$ 
and $E_{238}$, and hence having an $E$-irregularity index of 2)
(see \cite{em1}).
In 1954 L. Carlitz \cite{ca} proved that there are infinitely many 
$E$-irregular primes $p$, i.e., $p\mid E_2E_4\cdots E_{p-3}$. 
Using modular arithmetic to determine 
divisibility properties of the corresponding Euler numbers,
the $E$-irregular primes less than 10000 were found in 1978 by 
R. Ernvall and T. Mets\"{a}nkyl\"{a}  \cite{em1}.

In his book \cite[p. 203]{r} P. Ribenboim noticed that ``it is not 
all surprizing that the connection, via Kummer's theorem,
between the primes dividing certain Bernoulli numbers 
and the truth of Fermat's theorem, would suggest a similar
theorem using the Euler numbers." 
Vandiver \cite{va} proved that $x^p+y^p=z^p$ has no solution for 
integers $x,y,z$ with $\gcd(xyz,p)=1$ if $E_{p-3}\equiv 0\,(\bmod\,p)$.
The analogous result was proved by Cauchy (1847) and 
Genocchi (1852) (see \cite[p. 29, Lecture II, Section 2]{r}) with the 
Bernoulli number  $B_{p-3}$ instead of $E_{p-3}$.
Further, in 1950 M. Gut \cite{gu} proved that the condition $E_{p-3}\equiv 
E_{p-5}\equiv E_{p-7}\equiv E_{p-9}\equiv E_{p-11}\equiv 0\,(\bmod\,p)$ is 
necessary for the Diophantine equation $x^{2p}+y^{2p}=z^{2p}$ to be solvable. 

Furthermore, numerous combinatorial  congruences
recently obtained by Z.-W. Sun in \cite{s4}--\cite{s7} 
and by Z.-H. Sun in \cite{s2} 
involve  Euler numbers $E_{p-3}$ with a prime $p$. 
Many of these congruences become ``supercongruences"
if and only if $E_{p-3}\equiv 0\,(\bmod\,p)$ (A {\it supercongruence}
is a congruence whose modulus is a prime power.)
This gives a significance to primes $p$ for which 
$E_{p-3}\equiv 0\,(\bmod\,p)$. The first two primes 149 and 241 have 
also been discoverded by Z.-W. Sun \cite{s4}.

In this note, we focus  our attention to the computational search  
for residues of Euler numbers $E_{p-3}$ modulo a prime $p$. By the congruence
obtained in 1938 by E. Lehmer \cite[p. 359]{l}, for each prime $p\ge 5$
  \begin{equation}\label{con1}
\sum_{k=1}^{\left[p/4\right]}\frac{1}{k^2}\equiv (-1)^{(p-1)/2}4E_{p-3}
\pmod{p},
  \end{equation}
where $[a]$ denotes the integer part of a real number $a$.
Usually  (cf. \cite{em1}), if $E_{p-3}\equiv 0\,(\bmod\,p)$ then we say that 
$(p, p-3)$ is an $E$-{\it irregular pair}. It was founded in \cite{em1} that 
in the range $p<10^4$ $(p,p-3)$ is an $E$-irregular pair 
for $p=149$ and $p=241$.

For our computations presented in Section 3 we do not use 
Lehmer's congruence (1) including  harmonic number of the second order. Our 
computation via {\tt Mathematica 8} which uses the 
expression including the harmonic number  (of the first order)
is very  much faster than those related to the congruence (1). 
Here we report  that only  three primes less than $10^7$ 
satisfy the condition $E_{p-3}\equiv 0\,(\bmod\,p)$.
Using our computational results and statistical considerations similar 
to those  in relation to a search  for Wieferich and Fibonacci-Wieferich 
and Wolstenholme primes (cf. \cite[p. 447]{cdp} and \cite{mr}), 
we conjecture that there are infinitely many primes $p$ such that 
$E_{p-3}\equiv 0\,(\bmod\,p)$.

\section{A congruence used in our computation}

Here, as usually in the sequel, for integers $m,n,rs$ with
$n\not= 0$ and $s\not=0$, and a prime power $p^a$
we put $m/n\equiv r/s \,(\bmod{\,\,p^e})$ 
if and only if $ms\equiv nr\,(\bmod{\,\,p^e})$, and the residue
class of $m/n$ is the residue class of $mn'$ where 
$n'$ is the inverse of $n$ modulo $p^e$.  

In what follows $p$ always denotes a prime.
The Fermat Little Theorem states that if $p$ is a 
prime and $a$ is an integer not divisible by $p$,
then $a^{p-1}\equiv 1\,(\bmod{\,\,p})$. This gives rise
to the definition of the {\it Fermat quotient of $p$ to base $a$},
    $$
q_p(a):=\frac{a^{p-1}-1}{p},
    $$
which is an integer. 
It is well known that divisibility of Fermat quotient $q_p(a)$ by $p$
has numerous applications which include the Fermat Last Theorem and
squarefreeness testing (see \cite{em2}, \cite{gr1} and \cite{r}). 
If $q_p(2)$ is divisible by $p$, $p$ is said to be {\it Wieferich prime}. 
Despite several intensive searches, only two Wieferich primes 
are known: $p=1093$ and $p=3511$ (see \cite{cdp} and \cite{dk}).
Another class of primes initially defined because of Fermat Last Theorem
are {\it Fibonacci-Wieferich primes}, sometimes called {\it Wall-Sun-Sun 
primes}. A prime $p$ is said to be {\it Fibonacci-Wieferich prime} if 
the Fibonacci number $F_{p-\left(\frac{p}{5}\right)}$ is divisible by 
$p^2$, where $\left(\frac{p}{5}\right)$ denotes the Legendre symbol
(see \cite{sun}). 
A search in \cite{mr} and \cite{dk} shows that there are no 
Fibonacci-Wieferich primes less than $9.7\times 10^{14}$. 
 
For the computation of residues of Euler numbers $E_{p-3}$ 
modulo a prime $p$, it is suitable to use the following 
congruence which runs significantly faster than Lehmer's congruence (1).

\vspace{3mm} 
\noindent {\bf Theorem}
(\cite[Theorem 4.1(iii)]{s2}). {\it Let $p\ge 5$ be a prime. Then
  \begin{equation}\label{con4'}
\sum_{k=1}^{\left[p/4\right]}\frac{1}{k}+ 3q_p(2)
-\frac{3p}{2}q_p(2)^2\equiv (-1)^{(p+1)/2}pE_{p-3}\pmod{p^2},
    \end{equation}
where $[a]$ denotes the integer part of a real number $a$.} 

\begin{proof}
Quite recently, Z.-W. Sun 
\cite[Proof of Theorem 1.1, the congruence
after (2.3)]{s4} noticed that by a result of 
Z.-H. Sun \cite[Corollary 3.3]{s2}, 
   \begin{equation}\label{con4'}
\sum_{k=1}^{(p-1)/2}\frac{(-1)^{k-1}}{k}\equiv q_p(2)
-\frac{p}{2}q_p(2)^2-(-1)^{(p+1)/2}pE_{p-3}\pmod{p^2}.
    \end{equation}
On the other hand, we have
 \begin{equation}\label{con5'}
\sum_{k=1}^{(p-1)/2}\frac{(-1)^{k-1}}{k}=
\sum_{k=1}^{(p-1)/2}\frac{1}{k}-2\sum_{1\le k\le (p-1)/2\atop 2\mid k}
\frac{1}{k}=\sum_{k=1}^{(p-1)/2}\frac{1}{k}-\frac{1}{2}
\sum_{j=1}^{\left[p/4\right]}\frac{1}{j}.
   \end{equation}
By the classical congruence 
proved in 1938 by E. Lehmer \cite[the congruence (45), p. 358]{l},
for each prime $p\ge 5$ 
 \begin{equation}\label{con6}
\sum_{k=1}^{(p-1)/2}\frac{1}{k}\equiv -2q_p(2)+pq_p(2)^2\pmod{p^2}.
   \end{equation}
Substituting the congruence (5) into (4), we obtain 
 \begin{equation}\label{con5'}
\sum_{k=1}^{(p-1)/2}\frac{(-1)^{k-1}}{k}\equiv 
-2q_p(2)+pq_p(2)^2\pmod{p^2} -\frac{1}{2}\sum_{j=1}^{\left[p/4\right]}
\frac{1}{j}\pmod{p^2}.
   \end{equation}
Finally, substituting (6) into (3), we immediately obtain (2).
\end{proof}

\section{The computation}

Using the  congruence (2), a computation via {\tt Mathematica 8} shows that only  
three primes less than $10^7$ 
 satisfy the condition $E_{p-3}\equiv 0\,(\bmod\,p)$
(such primes are 149, 241 and 2946901, and 
they are given as a sequence A198245 in \cite{sl}).
Notice also that in 2011 \cite[p. 3, Remarks]{me}, the author
of this article reported that these three primes 
are only primes less than $3\times 10^6$. 

Recall that investigations of such primes 
have been recently suggested by Z.-W. Sun in \cite{s4}; namely, 
in \cite[Remark 1.1]{s4} Sun found the first and the second such primes,
149 and 241, and used them to discover curious 
supercongruences (1.2)--(1.5) from Theorem 1.1 in \cite{s4} involving 
$E_{p-3}$. 

Motivated by search for Wieferich and Fibonacci-Wieferich 
primes given in \cite{cdp} and \cite{dk}
and search for Wolstenholme primes given in \cite{mr}, 
here we use similar computational considerations for Euler numbers
$E_{p-3}$ where $p$ is a prime. Our computational results presented below
suggest  two conjectures on numbers $E_{p-3}$ that are analogous to those
on Wieferich (\cite{cdp}, \cite{dk}) and Wolstenholme primes \cite{mr}.
Accordingly, we search primes $p$ in the range $[10^5, 5\times 10^6]$
such that $E_{p-3}\equiv A\,(\bmod\,p)$ with $|A|\le 100$ and/or  
$10^4\cdot |A/p|\le 1$. Our search employed the congruence (2)
which runs significantly faster than Lehmer's congruence (1)
and than the code 
\vspace{2mm}

{\tt Print[\{Prime[n]\},Mod[EulerE[Prime[n]-3],Prime[n]]]}
\vspace{2mm}

Here {\tt EulerE[k]} gives $E_k$  and {\tt Mod[a,m]} gives $a(\bmod{\,m})$.

\noindent in {\tt Mathematica 8}, as well as than  some other known 
congruences involving Euler number $E_{p-3}$. 

Namely, in order to obtain data of Table 1 below 
concerning primes $p$ with  $10^5<p< 5\times 10^6$ we  used the code:
\vspace{2mm} 

\noindent{\tt{\small
 Do[If[Max[Min[Mod[(Mod[Numerator[HarmonicNumber[Floor[Prime[n]/4]]],

\noindent Prime[n]\^{}2]\*PowerMod[Denominator[HarmonicNumber
[Floor[Prime[n]/4]]],

\noindent -1,Prime[n]\^{}2+3*(2\^{}(Prime[n]-1)-1)/Prime[n]-PowerMod[2,-1,Prime[n]\^{}2]

\noindent *(3*Prime[n])*((2\^{}(Prime[n]-1)/Prime[n]\^{}2)/((-1)\^{}((Prime[n]+1)/2)

\noindent *Prime[n]),Prime[n]],Prime[n]-Mod[(Mod[Numerator[HarmonicNumber

\noindent[Floor[Prime[n]/4]]],Prime[n]\^{}]*PowerMod[Denominator[HarmonicNumber

\noindent[Floor[Prime[n]/4]]],-1,Prime[n]\^{}2+3*(2\^{}(Prime[n]-1)-1)/Prime[n]

\noindent -PowerMod[2,-1,Prime[n]\^{}2]*(3*Prime[n])*((2\^{}(Prime[n]-1)-1)
/Prime[n])\^{}2)

\noindent /((-1)\^{}((Prime[n]+1)/2)*Prime[n]),Prime[n]]]]==1000, Print[\{n\},

\noindent \{Prime[n]\},\{Mod[(Mod[Numerator[HarmonicNumber[Floor[Prime[n]/4]]],

\noindent Prime[n]\^{}2]*PowerMod[Denominator[HarmonicNumber[Floor[Prime[n]/4]]],-1,

\noindent Prime[n]\^{}2+3*(2\^{}(Prime[n]-1)-1)/Prime[n]-PowerMod[2,-1,Prime[n]\^{}2]

\noindent *(3*Prime[n])*((2\^{}(Prime[n]-1)/Prime[n]\^{}2)/((-1)\^{}((Prime[n]+1)/2)

\noindent *Prime[n]),Prime[n]]\},\{Prime[n]-Mod[(Mod[Numerator[HarmonicNumber

\noindent[Floor[Prime[n]/4]]],Prime[n]\^{}2]*PowerMod[Denominator[HarmonicNumber

\noindent[Floor[Prime[n]/4]]],-1,Prime[n]\^{}2+3*(2\^{}(Prime[n]-1)-1)/Prime[n]

\noindent-PowerMod[2,-1,Prime[n]\^{}2]*(3*Prime[n])*((2\^{}(Prime[n]-1)/Prime[n]\^{}2)

\noindent /((-1)\^{}((Prime[n]+1)/2)*Prime[n]),Prime[n]]\}]],\{n,i,j\}]}}

\vspace{2mm}

Here {\tt Mod[a,m]} gives $a(\bmod{\,m})$, 
{\tt PowerMod[a,b,m]} gives $a^b(\bmod{\,m})$ (and is faster than
{\tt Mod[a\^{}b,m]}.

Further, in order to verify that there are no primes $p$ between 
$5\times 10^6$ and $10^7$ such that $E_{p-3}\equiv 0\,(\bmod\,p)$,
 we used  the following code which is very  much faster the  previous code:
\vspace{2mm} 
  
\noindent {\tt\small Do[Print[\{n\},\{Prime[n]\},Mod[Numerator[2*HarmonicNumber[Floor[

\noindent Prime[n]/4]]+6*(2\^{}(Prime[n]-1)-1)/Prime[n]-3*(2\^{}(Prime[n]-1)-1)\^{}2

\noindent/Prime[n]],Prime[n]\^{}2]],\{n,i,j\}]}

\vspace{2mm}
 
Certainly $A=A(p)$ can take  any of $p$ values $(\bmod{\,p})$.
Assuming that $A$ takes these values these values randomly, 
the ``probability" that $A$ takes any particular value
(say 0) is $1/p$. From this, in accordance to the 
heuristic given in \cite{cdp} related to the Wieferich primes,
we might argue that  the  number of primes $p$ 
in an interval $[x,y]$ such that $E_{p-3}\equiv 0\,(\bmod\,p)$  is
 expected to be
       \begin{equation}\label{con7}
\sum_{x\le p\le y}\frac{1}{p}\approx \log\frac{\log y}{\log x}.
      \end{equation}
If this is the case, we would be only expect to find
about $0.998529(\approx 1)$,  such primes in  the interval $[10^7,10^{19}]$.
On the other hand, since 9999991 is the greatest prime less than 
$10^7$ and is is actually 664589th prime,  by above estimate, 
we find that in the interval $[2,10^{7}]$ we can expect about 
$\sum_{2\le p\le 10^7}1/p=\sum_{k=1}^{664589}1/p_k\approx 3.04145$
primes $p$ such that $E_{p-3}\equiv 0\,(\bmod\,p)$ ($p_k$ 
is a $k$th prime); as noticed  previously, 
our computation shows that all these primes are 149, 241 and 2946901.
  \vfill\eject

\begin{center}{{\bf Table 1.}} 
Primes $p$ with  $10^5<p<5\times 10^6$ for which $E_{p-3}\equiv A\,(\bmod\,p)$
with $|A|\le 100$ and/or 
with related values $|A/p|\le 10^{-4}$  (given in multiples of $10^{-4}$)
 \end{center}
  \begin{tabular}{ccc}    
$p$ & $A$& $|A/p|$  \\\hline
105829 &- 74 & $>1$ \\
111733 & 45  & $>1$ \\
127487 & 38   & $>1$ \\
130489 & -27     & $>1$  \\
131617 & 9    &  0.683802\\
162847 & -85      &  $>1$\\
165157 &-46       &  $>1$\\
171091 &- 17      &  0.993623\\
171449 & 7       &  0.408285\\
191237 &37       &  $>1$\\
192961 & 63      &  $>1$\\
200461 &7       &  0.349195\\
209393 &27      &  $>1$\\
245471 &39      &  $>1$\\
246899 &-54     &  $>1$\\
276371 &-69     &  $>1$\\
290347 & 10     &  0.344415\\
292183 & 53     &  $>1$\\
306739 &-42     &  $>1$\\
317263 &-35     &  $>1$\\
321509 &84      &  $>1$\\
342569   &25    &  0.729780\\
422789   &-40    &0.946098  \\
429397   &-62  &  $>1$    \\
440047 & 82 & $>1$\\
479561 & 31 &0.646425\\
501317 & 60 & $>1$\\
546631 &92 & $>1$\\
628301 & 73& $>1$\\
636137 &25 & 0.392997\\
656147& -68& $>1$\\
659171 & -22 & 0.333753\\
687403 &-4 & 0.058190\\
717667 &-42 &0.585230\\
719947 &53 &0.736165\\
766261 &-8 &0.104403\\
801709 &53 & 0.661088\\
920921& -82 & 0.890413\\
924727 & -8 &   0.086512 \\
1064477 & $106(>100) $ &  0.995794 \\
1080091 & 42 & 0.388856\\
1159339 & -38 & 0.327773\\ 
1202843 & 21 &0.174586\\
1228691 & 15 &0.122081\\
1285301 &47 & 0.365673 \\
1336469 & -5 & 0.037412\\
1353281 & 78 &  0.576377\\
   \end{tabular}
 \begin{tabular}{ccc}    
$p$ & $A$& $|A/p|$  \\\hline
1355269 & -60 & 0.442717\\
1392323 & -29 &0.208285\\
1462421 & -78 & 0.533362\\
1546967 & -43 &   0.277963\\
1743271 & $107(>100)$ & 0.613789\\
1794049 & $-131(<-100)$  & 0.730192\\
1808497 & $-121(<-100)$  & 0.669109\\
1952131 & $-153(<-100)$   & 0.783759\\
1986539 & $-157(<-100)$   & 0.790319\\
2053873 & 18 & 0.087639\\
2114251 & $211(>100) $ & 0.997989\\
2236349 & 4 & 0.017886\\
2342381 & $143(>100)$ &  0.610490\\ 
2410627 & $-219(<-100)$   & 0.908477\\
2472731 & $230(>100)$ & 0.930146\\
2583011 & $159(>100)$  & 0.615561\\
2619847 & $224(>100)$  & 0.855011\\
2740421 & 225 & 0.821042\\
2890127 &-34 & 0.117642\\
{\bf 2946901} &{\bf 0} & {\bf 0} \\
3279833 &$-111(<-100)$   &  0.338432\\
3290689 & $200(>100)$  & 0.607775\\
3312653 & $228(>100)$  & 0.688270\\
3340277 & $226(>100)$  & 0.676591\\
3355813 & $116(>100)$  & 0.345669\\
3652613 & $-290(<-100)$  & 0.793952\\ 
3818131 & $-318(<-100)$  & 0.832868\\
3852677 & 75 & 0.194670\\ 
3960377 & -48 & 0.121201\\ 
4007747 & $190(>100)$  & 0.474082\\
4121503 &$-270(<-100)$   &0.655101  \\
4171229 & $153(>100)$  & 0.366798\\
4343659 &$-252(<-100)$   &  0.580156\\
4392007 & 55 &  0.125227 \\
4418497 & 70 & 0.158425\\
4475707 & $193(>100)$ & 0.431217\\
4541501 & $120(>100)$ & 0.264230\\
4551973 & $-362(-100)$  & 0.795260\\
4564939 & -63 & 0.138008\\
4631399 & $367(>100)$ & 0.792417\\
4674347 & $302(>100)$ & 0.646080\\
4706047 & $220(>100)$ & 0.467484\\
4751599 & $-279(<-100)$  & 0.587171\\
4761677 & $200(>100)$ & 0.420020\\
4869517 & -100 & 0.205359\\
4898099 & $-236(<-100)$  & 0.481820\\
4928503 & $-173(<-100)$  & 0.351019\\
     \end{tabular}}
\vfill\eject

The second column of Table 1 shows that there are 
61 primes between $10^5$ and $5\times 10^6$ for which $|A|\le 100$.
Since the ``probability" that $|A|\le 100$ for a prime $p\gg 200$ 
is equal to    $201/p$,
it follows that expected number of such primes between 
$M$th prime $p_M$ and  $N$th prime $p_N$ with $N>M\gg 1000$  (that is, 
$p_N>p_M\gg 1000$) is equal to
  \begin{equation}\label{con8}
Q(N,M,100)=201\sum_{p_M< p<p_N}\frac{1}{p},
 \end{equation}
 where the summation ranges over all primes $p$ such that 
$p_M< p<p_N$. In particular, for the values $M=9593$ and $N=348513$ 
which correspond to the interval $[10^5,5\times 10^6]$ 
containing  all primes from Table 1, we have
  \begin{equation}\label{con9}
Q(348513,9593,100)=
201\sum_{10^5< p<5\times 10^6}\frac{1}{p}\approx
201\cdot 0.292251=58.742451.
 \end{equation}

On the other hand, Table 1 shows that there are 
61 primes between $10^5$ and $5\times 10^6$ for which $|A|\le 100$,
which is $\approx 3.8431\%$ greater than  related 
``expected number" $58.742451$.

Because our program recorder all $p$ with ``small $|A|$", that is,
with $|A|\le 100$, we compiled a large data set which can be used to give
more rigorous (experimental) confirmation of both our Conjectures 1
and 2. Indeed, our program recorded 568  primes $p$ in 
the interval $[10^5,5\times 10^6]$  for which 
$|A|\le 1000$. On the other hand, according to the 
formula (9), it follows that expected number of such primes 
is equal to
 \begin{equation}\label{con10}
Q(348513,9593,1000)=
2001\sum_{10^5< p<5\times 10^6}\frac{1}{p}\approx
2001\cdot 0.292251=584.794251
 \end{equation} 
which is $\approx 2.956 \%$ greater than related ``expected number" 568. 

Instead, of selecting values based on $|A|\le 100$,
we suggest to select them based on $A/p <q\times 10^{-4}$
(e.g., $q=1$) that would be consistent with the original selection criterion.
In particular, in the third column of Table 1
there are 72 primes $p$ contained in the interval  $[10^5,5\times 10^6$ 
with related values $10^4\times A/p <1$.

Furthermore, since the ``probability" that $|A/p|\le 10^{-4}$ for a prime 
$p\gg 10000$ is equal to 
  $$
\frac{2\left[\frac{p}{10000}\right]+1}{p}\approx \frac{2}{10000},
 $$
it follows that expected number of such primes between 
$M$th prime $p_M$ and  $N$th prime $p_N$ with $N>M\gg 1000$  (that is, 
$p_N>p_M\gg 10000$) is equal to
  $$
P(N,M)=\frac{2(N-M)}{10000}.
 $$
In particular, for the values $M=9593$ and $N=348513$ 
which correspond to the range $(10^5,5\times 10^6)$ 
of all primes from Table 1, we have
  $$
P(348513,9593)=\frac{677840}{10000}=67.7840 
 $$
which is $\approx 5.855\%$ less than 72.

All the previous considerations and the well known fact that the series 
 $$
\sum_{p\,\,\rm{ prime}}\frac{1}{p}
 $$
diverges suggest the following conjecture.
\vspace{2mm}

\noindent {\bf Conjecture 1.}  There are infinitely many primes 
$p$ such that $E_{p-3}\equiv 0\,(\bmod\,p)$.
\vspace{2mm}

Since 
   $$
\sum_{x\le  p\le y}\frac{1}{p}\approx \log\log x-\log\log y,
   $$
in view of the previous comparison of our computational results 
with expected number of  primes $p\in [10^5,5\times 10^6]$  for which 
$|A(p)|\le 100$ given by (9) (or primes $p\in [10^5,5\times 10^6]$  for which 
$|A(p)|\le 1000$ given by (10)),  we can assume that
expected number of  primes $p$ in an interval $[x,y]$ 
such that $K\le |A(p)|\le L$ is asymptotically equal to 
(cf. (7))
  \begin{equation}\label{con11}
2(L-K)\cdot (\log\log b-\log\log a).
 \end{equation}
Using a larger data set which our program recorded, 
consisting of total  568  pairs $(p,A(p))$ such that 
$p\in [10^5,5\times 10^6]$ and $|A(p)|\le 1000$, we obtain
 experimental results presented in Table 2. In Table 2 the values
in ``column $k$" and  in first and second 
row  reflect the number of $p\in [10^5,10^6]$
and $p\in [10^6,5\times 10^6]$, respectively,  
such that $A=A(p)\in [k\times 100,(k+1)\times 100]$ ($k=0,1,\ldots 9$). 
Expected numbers given in the last column of Table 2 are calculated by the 
formula (11).  

\begin{center}{{\bf Table 2.}}
 \end{center}
\vspace{-1mm}
  \begin{tabular}{|c|cccccccccc|c|}\hline    
 &&&&& $k$ &&&&&& \\\hline
  Interval & 0& 1& 2& 3& 4& 5& 6& 7 & 8 & 9& Expected\\\hline                
$[10^5,10^6]$ &42 &51 &37 &30 &29 &24 &31 &34  &42  &44 & 36.464\\
$[10^6,5\times 10^6]$ & 22& 23& 26& 20& 22& 22& 21&  24&  21& 20& 22.039
\\\hline
 \end{tabular}\\

 Table 2 presents a small snapshot of our experimental results. 
Notice that by the data of the last row, the relative error between 
the conjectured and experimental values for $k=0,1,\ldots ,9$
are respectively equal to $0.18\%$, $4.18\%$, $15.23\%$, 
$10.20\%$, $0.18\%$, $0.18\%$, $4.95\%$, $8.17\%$,
$4.95\%$, $10.20\%$. 
Accordingly, we propose the following conjecture (cf. the same conjecture
in \cite[Conjecture 6.1]{dk} concerning the Wieferich primes; 
see also \cite[Section 3]{cdp}).
 
\vspace{2mm}
\noindent {\bf Conjecture 2.} The number of primes 
$p\in [a,b]$ such that $|A|=|A(p)|\in [K,L]$ is {\it asymptotically}
\vspace{2mm}
  $$
2(L-K)\cdot (\log\log b-\log\log a).
 $$

\vspace{2mm}
\noindent {\bf Remarks.}
Recall that a prime $p$ is said to be a {\it Wolstenholme prime} if it 
satisfies the congruence 
 $$
{2p-1\choose p-1} \equiv 1 \pmod{p^4},
 $$
or equivalently (cf. \cite[Corollary on page 386]{mc}; also 
see \cite{gl}) that $p$ divides the numerator of $B_{p-3}$.
The only two known such primes are 16843 and 2124679, and 
by a result of R.J. McIntosh and E.L. Roettger from 
\cite[pp. 2092--2093]{mr}, these primes are the only two  Wolstenholme 
primes less than $10^9$.
Nevertheless, using similar arguments to those given in Section 3 of this 
paper, McIntosh  \cite[page 387]{mc} conjectured that there are infinitely many 
Wolstenholme primes. 

\vfill\eject


{\bf\large References}

\begin{thebibliography}{99} 

\bibitem{ca} L. Carlitz, Note on irregular primes,  {\it Proc. Amer. Math. Soc.}
{\bf 5} (1954) 329--331.

\bibitem{cdp} R. Crandall, K. Dilcher and C. Pomerance,
A search for Wieferich and Wilson primes,
{\it Math. Comp.} {\bf 66} (1997) 443--449. 

\bibitem{dk} F.G. Dorais and D. Klyve,  A Wieferich prime search up to
$6.7\times 10^{15}$,  {\it J. Integer Seq.} {\bf 14} (2011) Article 11.9.2.
 
\bibitem{em1}  R. Ernvall and T. Mets\"{a}nkyl\"{a}, Cyclotomic invariants and
 $E$-irregular primes,  {\it Math. Comp.} {\bf 32} (1978) 617--629. 

\bibitem{em2}  R. Ernvall and T. Mets\"{a}nkyl\"{a}, On the 
$p$-divisibylity of Fermat quotients, 
{\it Math. Comp.} {\bf 66} (1997) 1353--1365. 

\bibitem{gl} J.W.L. Glaisher,  On the residues of the sums of 
products of the first $p-1$ numbers, and their powers, to modulus $p^2$
or $p^3$,   {\it Q.  J. Math.} {\bf 31} (1900) 321--353. 

\bibitem{gr1} A. Granville,   {\it Some conjectures related to 
Fermat's Last Theorem, Number Theory} (Banff, AB, 1988),
 de Gruyter, Berlin, 1990, 177--192.

\bibitem{gu} M. Gut, Eulersche Zahlen und grosser Fermat'scher Satz,
{\it Comment. Math. Helv.} {\bf 24} (1950) 73--99.


\bibitem{l} E. Lehmer, On congruences 
involving Bernoulli numbers and the quotients of
Fermat and Wilson,   {\it Ann. Math.} {\bf 39} (1938) 350--360. 



\bibitem{mc} R.J. McIntosh,  On the converse of 
 Wolstenholme's theorem,  {\it Acta Arith.}  {\bf 71} (1995) 381--389. 

\bibitem{mr} R.J. McIntosh and E.L. Roettger, A search for 
Fibonacci-Wieferich and Wolstenholme primes, {\it Math. Comp.}
{\bf 76} (2007) 2087--2094.

\bibitem{me} R. Me\v{s}trovi\'{c}, An exstension of a congruence by 
Kohnen, 13 pages, preprint {\tt arXiv:1109.2340v3 [math.NT]} (2011). 

\bibitem{r} P. Ribenboim, 13   {\it Lectures on Fermat's Last Theorem}, 
Springer-Verlag, New York, Heidelberg, Berlin, 1979.

\bibitem{sl} N.J.A. Sloane, {\it Sequence} A198245  
in OEIS (On-Line Encyclopedia of Integer Sequences), 
{\tt http://oeis.org/A198245}.


\bibitem{sc} H.M. Srivastava and J. Choi,  {\it Series Associated with the
Zeta and Related Functions}, Kluwer Academic Publishers, Dordrecht, Boston
and London, 2001. 


\bibitem{s1} Z.-H. Sun, Congruences concerning Bernoulli numbers
and Bernoulli polynomials,  {\it Discrete Appl. Math.} {\bf 105} (2000)
193--223.

\bibitem{s2} Z.-H. Sun, Congruences involving Bernoulli and Euler numbers, 
{\it J. Number Theory}, {\bf 128} (2008) 280--312.

\bibitem{sun} Z.-H. Sun and Z.-W. Sun, Fibonacci numbers and Fermat's last
theorem,   {\it Acta Arith.} {\bf 60} (1992) 371--388.

\bibitem{s3} Z.-W. Sun, Binomial coefficients, Catalan numbers 
and Lucas quotients,  {\it Sci. China Math.} {\bf 53} (2010) 2473--2488;
preprint {\tt arXiv:0909.5648v11 [math.NT]} (2010).

\bibitem{s4} Z.-W. Sun, On Delannoy numbers and Schr\"{o}der numbers,
{\it J. Number Theory} {\bf 131} (2011) 2387--2397;
preprint {\tt arXiv:1009.2486v4 [math.NT]} (2011).

\bibitem{s5} Z.-W. Sun, Super congruences and Euler numbers,
{\it Sci. China Math.} {\bf 54} (2011) 2509--2535; 
preprint {\tt arXiv:1001.4453v19 [math.NT]} (2011). 

\bibitem{s7} Z.-W. Sun, On congruences related to central 
binomial coefficients,  {\it J. Number Theory} {\bf 131} (2011) 2219--2238;
preprint {\tt arXiv:0911.2415v16 [math.NT]} (2011).

\bibitem{s6} Z.-W. Sun, A refinement of a congruence result 
by van Hamme and Mortenson, accepted for publication in 
{\it Illinois J. Math.}; preprint {\tt arXiv:1011.1902v5 [math.NT]} (2011). 

\bibitem{va} H.S. Vandiver, Note on Euler number criteria for the
first case of Fermat's last theorem,  {\it Amer. J. Math.} {\bf 62}
(1940) 79--82.
\end{thebibliography}
\end{document}